\documentclass[a4paper,12pt]{article}

\usepackage[cp1251]{inputenc}
\usepackage[english]{babel}

\usepackage{amssymb,amsmath,amsthm,cite}
\usepackage[pdftex]{hyperref}
\usepackage{graphics}

\usepackage[centerlast,small]{caption2}
\begin{document}
\def\leq{\leqslant}
\def\le{\leqslant}
\def\geq{\geqslant}
\def\ge{\geqslant}

\renewcommand{\thefigure}{\arabic{figure}}
\newcommand{\Cb}{\mathbb{C}}

\newtheorem{te}{Theorem}
\newtheorem{pr}{Proposition}
\newtheorem{lem}{Lemma}
\newtheorem{sle}{Corollary}
\theoremstyle{definition}
\newtheorem{za}{Remark}
\newtheorem{opr}{Definition}
\newtheorem{ex}{Example}


\def\inserthyphen{\ifcat\next a-\fi\ignorespaces}
\def\pblack-{$\bullet$\futurelet\next\inserthyphen}
\def\pwhite-{$\circ$\futurelet\next\inserthyphen}
\def\pcross-{$\vcenter{\hbox{$\scriptstyle\times$}}$\futurelet\next\inserthyphen}
\def\black{\protect\pblack}
\def\white{\protect\pwhite}
\def\cross{\protect\pcross}

\newcommand{\Cyl}{\operatorname{Cyl}}
\newcommand{\Int}{\operatorname{Int}}
\newcommand{\ind}{\operatorname{ind}}
\newcommand{\prr}{\operatorname{pr}}

\title{A cellular structure of the space of branched coverings of the two-dimensional sphere}
\author{V.~I.~Zvonilov 
\thanks{The work of the first author was done on a subject of the State assignment N 0729-2020-0055}, S.~Yu.~Orevkov}



\date{January 11, 2020}

	
\maketitle

\begin{abstract}
 
 For a closed oriented surface $ \Sigma $, let $ X _ {\Sigma, n} $ be the space of isomorphism classes of $ n $-fold orientation preserving branched coverings $ \Sigma \rightarrow S^2 $ of the two-dimensional sphere. Earlier, the authors constructed a compactification $ \bar {X} _ {\Sigma, n} $ of this space, which coincides with the Diaz-Edidin-Natanzon-Turaev compactification of the Hurwitz space $ H (\Sigma, n) \subset X _ {\Sigma , n} $ that consists 
  of isomorphism classes of branched coverings with all critical values being simple.
 With the help of Grothendieck's {\it dessins d'enfants}, a cellular structure of this compactification is constructed. The results obtained are applied to the space of trigonal curves on an arbitrary Hirzebruch surface.
 
\end{abstract}

\section {Introduction}
Let $ \Sigma $ be a closed oriented surface (fixed throughout the article). Orientation preserving $ n $-fold branched coverings $ f_1, f_2 \colon \Sigma \rightarrow S^2 $ of the two-dimensional sphere are called \emph{isomorphic} if there is a
homeomorphism $ \alpha \colon \Sigma \rightarrow \Sigma $ with $ f_2 \circ \alpha = f_1 $.
Let $ X _ {\Sigma, n} $ be the set of isomorphism classes of such coverings. The paper \cite {ZO} introduces a topology on this set and constructs a compactification $ \bar {X} _ {\Sigma, n} $ of the resulting space, which coincides with the Diaz-Edidin-Natanzon-Turaev compactification of the Hurwitz space $ H ( \Sigma, n) \subset X _ {\Sigma, n} $ consisting of isomorphism classes of coverings with simple critical values. A point of the space $ \bar {X} _ {\Sigma, n} $ is the  isomorphism class of a degeneration $ f '\colon \Sigma' \rightarrow S^2 $ (see section \ref {Deform}) of a branched covering $ f \colon \Sigma \rightarrow S^2 $, where $ \Sigma '$ is a degeneration of the surface $ \Sigma $ (see Section \ref {singsurf}).

The main result of this work is the introduction of a cellular structure, i.e.  the structure of the $ CW $-complex, into $ \bar {X} _ {\Sigma, n} $, using the notion of graph of branched covering  of the sphere.
In the constructed cellular space, cells are the isomorphism classes of branched coverings of the sphere with isomorphic graphs.
To construct the graph of branched covering  of the sphere, we follow the principle formulated in \cite [Principle 1.6.1] {LZ}, which proposes to take as  the graph on the covering surface the preimage of a base graph on the sphere containing all critical values of the covering. The choice of the base graph is determined by the further application of the cellular structure.

In \cite {DE}, a cellular decomposition of the space $ \bar {X} _ {\Sigma, n} $ was constructed, which was used by the authors to calculate the homology of the space $ H (\Sigma, n) $. The authors did not prove that their partition is a $ CW $-complex (see \cite [Section 4.4] {DE}). However, replacing the base graph (see Remark \ref {basegraph} below) allows us to prove that this is true, in the same way we do it for our decomposition.

We apply the results obtained to the compactification of the space of $ j $-invariants of trigonal curves on a ruled
surface, i.e.  mappings of the base of the ruled surface
into a modular curve, see \cite [Section 4] {Or}, and also \cite [2.1.2, 3.1.1] {Degt}.
For working with trigonal curves, we choose a base graph on the sphere, which we call a topological hosohedron (see Subsection \ref {Grf}). The cellular decomposition of the work \cite {DE}, which is more economical in the number of cells than ours, is not convenient for studying trigonal curves, since the $ j $-invariant has fixed critical values.

The results obtained allow us to propose a way for computing the
fundamental group of
the space of nonsingular trigonal curves; knowledge of the group will help to obtain new restrictions on the topology of real algebraic varieties, in particular, of surfaces of degree 5 in real projective space.

The structure of the paper is as follows. Section \ref{space} contains the necessary information from \cite {ZO}: in Subsections \ref {singsurf}, \ref {Deform}  the concepts of a degenerate (singular) surface 
and its covering of the two-dimensional sphere are reminiscent;  in Subsection \ref {Xn} we repeat the definition of topology in $ \bar {X} _ {\Sigma, n} $.   In Subsections \ref {Grf}, \ref {Transform} graphs of branched coverings and transformations of the graphs are studied. In Section \ref {Cell} a cellular structure is introduced and a dual partition of the space $ \bar {X} _ {\Sigma, n} $ and its subspaces is constructed. In Section \ref {Hirz} the results obtained are applied to the space of trigonal curves.
\vspace{-1mm}
\section {The space of isomorphism classes of branched coverings.} \label{space}
Let's repeat the necessary
information from the paper \cite {ZO}.
\nopagebreak
\subsection {Singular surface}
\label{singsurf}

Everywhere below,  \emph{a surface}  is a Hausdorff topological space with a countable base such that any its point has a neighbourhood homeomorphic to an open disk  or to a wedge of open disks (or to the union of an open half-disk  with its diameter in case of a surface with boundary).  We call such a neighbourhood \emph{admissible} and say that the disks of the wedge are \emph{the branches of the surface} at the center of the wedge.
We assume that a non-negative integer 
$ g_v$ (called  \emph{the local genus of the surface at the point} $v$) is assigned to any point $v$ of such a surface; the local
genus is non-zero only at a finite number of points and it vanishes on the boundary.  That is,  speaking more formally, a surface is a pair of a topological space  and an integral-valued function  $v\mapsto g_v$ on it such that the space and the function have the above properties.  
Let $m_v$ be the number of branches of a surface at the point  $ v $ and  $ \mu_v=2g_v+m_v-1$ be  \emph{the Milnor number} of $ v $.
Points with  
$ \mu_v>0$ are \emph{the singular points} of the surface. A surface  is 
\emph{degenerate} (or \emph{singular}) if is has
singular points. Otherwise it is  \emph{nonsingular}.

The \emph{the normalization} of a surface is defined as a nonsingular
surface obtained by replacing an admissible neighbourhood of each
singular point with a disjoint union of smooth disks. There is a natural projection of the normalization to the surface that takes each glued disk  to the corresponding disk  of the wedge. 
\emph{A component} of a surface is a connected component of its normalization. So    \emph{a component of a surface is a two-dimensional manifold.} A surface is \emph{oriented/closed} if its normalization is oriented/closed (compact and without  boundary).

Below all the surfaces are oriented and all the maps are
orientation preserving. 

Let $ \Sigma_0 $ be a surface (in the above sense, maybe singular), $ v $ be its point, and
$U_v$ be the closure of an admissible neighbourhood of $v $.
Take a compact connected orientable  (maybe singular) surface 
$\tilde{U}_v$ with $m_v$ boundary components and such that
$b_1(\tilde{U}_v)+\sum_{x\in\tilde{U}_v}\mu_x=\mu_v$ and $b_2(\tilde{U}_v)=0$ (i.e. $\tilde{U}_v$ has no components without boundary);
here $b_i$ is $i$-th Betti number. Glue  the surfaces $ \Sigma_0\setminus\mathrm{Int}U_v $ and $\tilde{U}_v$ along their common boundary. We say that the resulting surface 
$ \Sigma_1$ is \emph{a perturbation} of $ \Sigma_0$ \emph{at the point} $ v $, or \emph{a local perturbation}, and $ \Sigma_0$ is \emph{a local degeneration} of $ \Sigma_1$, \emph{by means of
} $U_v$, $\tilde{U}_v$.
A composition of a finite number of local perturbations/degenerations of a surface is called a \emph{perturbation/degeneration} of the surface.

 \subsection{Perturbations/degenerations of coverings}
\label{Deform}
A mapping $ f $  of a surface $ \Sigma' $ to the $ 2 $-sphere is called  \emph{a branched covering} if, for $ \pi $ being the projection of normalization, the restriction of the composition $f\circ\pi$  to any component of the surface is an orientation preserving branched covering and each singular point   $ v \in\Sigma' $ with $ m_v=1 $  
is a ramification point of 
$f$. The ramification points of $f$ and the singular points of
$\Sigma'$ are called  \emph{critical points} of $f$ and their images are called
\emph{critical values}.  By $ \mathrm{cr} f $ we denote the set of all critical values of $ f $.

Let $v$ be a critical point of a branched covering $f_0:\Sigma_0\rightarrow S^2$
and let $\bar{U}_w\subset S^2$ be a closed disk which contains
$w=f_0(v)$ as an inner point and which does not contain other critical
values of $f_0$. 
Pick a connected component  $U_v$  of $f^{-1}(\bar{U}_w)$ with $v\in U_v$. Clearly it is a closed admissible neighbourhood of  $v$.
Let $ \Sigma_1 $ be a perturbation of $ \Sigma_0 $  at $ v $ by means of
$U_v$, $\tilde{U}_v$. We say that a branched covering  $f_1:\Sigma_1\rightarrow S^2$ is  \emph{a perturbation} of $f_0$ \emph{at} $ v $, or \emph{local perturbation with  perturbation domain} $ \bar{U}_w $ and that  $f_0$ is  \emph{a local degeneration} of $f_1$ if  $f_0=f_1$ on $\Sigma_0\setminus  \mathrm{Int}U_v=\Sigma_1\setminus\mathrm{Int}\tilde{U}_v$. It is obvious that $ f_1(\tilde{U}_v)=\bar{U}_w $.

Now suppose that the disks $\bar{U}_w$ picked for all critical values $ w $ of  $ f_0$ are disjoint.  
A composition $ f$ of a finite number of local perturbations of $ f_0$ that have the perturbation domains $ \bar{U}_w $ is called   \emph{a perturbation} of the covering. The union $V=\bigcup_w\bar{U}_w$ is \emph{a perturbation domain} of $f_0$. In this case  $ f_0 $ is called \emph{a degeneration} of $ f $.

Below, we need degenerations only of the fixed   \emph{nonsingular} surface $ \Sigma$ 
(see Introduction). 

Let  $ f':\Sigma'\rightarrow S^2$ be a degeneration of a branched covering $ f:\Sigma\rightarrow S^2$ corresponding to a perturbation domain $ V $. For any critical point  
$ v$ of $ f'$, let $ U_v\subset\Sigma',\tilde{U_v}\subset\Sigma$ be the corresponding surfaces involved in the definition
of perturbation.  
Denote by $k_1,k_2,\ldots,k_{m_v}$ the degrees of the restrictions of $ f' $ to the boundary circles of   $U_v$. It is clear that  $\sum_{i=1}^{m_v}k_i$ is equal to the local degree $ \mathrm{deg}_vf' $ of $ f' $ at the point $ v $ that is the multiplicity of the point as a root of the equation $w=f'(v)$. 
\emph{The index of a point} $ v\in \Sigma' $ is defined as $\mathrm{ind}_v=\mu_v + \mathrm{deg}_vf'$. 

\emph{The multiplicity of a critical value} $ w $ of  $ f' $ is $ \sum_{v\in f'^{-1}(w)}(\mathrm{ind}_v-1) $. A critical value is called \emph{simple} if it is of multiplicity $1$. Otherwise it is called \emph{multiple}.

It is clear that if two coverings are isomorphic then for any perturbation/degeneration $ f $ of  one of them there is a perturbation/degeneration of another one isomorphic to $ f $. So we can speak about a perturbation/degeneration of an  isomorphism class  of covering.

Let $V=\bigcup_{e\in E} D_e\subset S^2$ be the union of disjoint closed topological disks (here $e$ is a fixed interior point of the disk $D_e$).  Suppose that all the critical values of branched coverings  
$ f_0:\Sigma_0\rightarrow S^2$ and $ f_1:\Sigma_1\rightarrow S^2$ are simple and lie in $\mathrm{Int}\, V $.
We say that these coverings are    $(V,E)$-\emph{equivalent} if there exist homeomorphisms    
$ \alpha:\Sigma_0\rightarrow \Sigma_1$ and $\varphi:S^2\rightarrow S^2$ such that:   
\renewcommand{\labelenumi}{(\theenumi)}
\renewcommand{\theenumi}{\roman{enumi}}
\begin{enumerate}
	\item $ \varphi\circ f_0=f_1\circ\alpha $;
	\item  $ \varphi $ leaves the points of  $E\cup~(S^2~\setminus~\mathrm{Int}V)$ fixed and it is isotopic to the identity mapping in the class of such homeomorphisms leaving the points fixed.
	
\end{enumerate}
\begin{pr}[{see \cite [Proposition 4] {ZO}}] \label{pert}
\hspace{-0.8mm}For any branched covering \hspace{-0.5mm}$ f_0\!:\!\Sigma_0\!\rightarrow\!S^2$ there exists its perturbation $f_1:\Sigma_1\rightarrow S^2$  at a critical point $ v $ such that the surface  $\tilde{U}_v\subset \Sigma_1$ (involved in the definition of perturbation)  is nonsingular and  $\chi(\tilde{U}_v)=2-2g_v-m_v$. Moreover, $ f_1$ can be chosen so that all of its critical values  belonging to the perturbation domain $ \bar{U}_w $, $ w=f_0(v) $, are simple and  can be placed at any prescribed positions;
in this case the covering $ f_1 $ is unique up to $ (\bar{U}_w,w) $-equivalence.\qed
\end{pr}

\subsection{Topology on the set of isomorphism classes of branched coverings}
\label{Xn}
Remind that $ \Sigma $ is a fixed nonsingular oriented closed surface. As in  Introduction, let  $X_{\Sigma,n}$ denote the set of isomorphism classes of $ n $-folded branched coverings $ \Sigma\rightarrow S^2 $ and $\bar{X}_{\Sigma,n} $ denote the set of all degenerations of isomorphism classes of coverings from $X_{\Sigma,n}$. 

Given  $f$ with $ [f]\in\bar{X}_{\Sigma,n} $ and a perturbation domain $ V $ of $ f $,  
let $ U_{[f],V}\subset\bar{X}_{\Sigma,n} $ be the set  of isomorphism classes of all perturbations of $ f $ whose perturbation domain is $ V $. In other words, these are the classes of perturbations of coverings in $ [f] $ such that all their critical values belong to $ V $.

\begin{pr} [{see \cite [Proposition~5] {ZO}}] \label{base}
	The family  $$
\{U_{[f],V}\mid \text{$f\in\bar{X}_{\Sigma,n}$ \emph{and $V$ is a perturbation domain of} $f$}  \}
$$
is a base of topology on  $\bar{X}_{\Sigma,n} $. The subset $X_{\Sigma,n}$ is dense in the space  $\bar{X}_{\Sigma,n} $ endowed with this topology.\qed
	\end{pr}

\section {Hosography}
\subsection {The graph of branched covering of the sphere}
\label{Grf}

Let $ a $ be a positive integer.
Choose two diametrically opposite points on the two-dimensional sphere $ S^2 $ named \emph{poles}, mark one pole in black, the other one as a cross and call them, respectively, \black--\emph{vertex} and \cross--\emph{vertex}. Connect these vertices by $ a $ semicircles named \emph{meridians}, orient the meridians in the direction from the \black-- to the \cross--vertex and mark them with the integers from $ 1 $ to $ a $ when bypassing the \black--vertex  counterclockwise. Let us choose a finite set of points on (not necessarily all) the meridians. We say that an oriented graph $ \Theta $ is a
 \emph{topological hosohedron}  or simply \emph{hosohedron} if its vertices are the
poles and  the selected points, its edges are the arcs of the meridians, oriented in the same way as the meridians, and each edge are marked with the number of the meridian containing the edge. The \black-- and \cross-- vertices are  \emph{stationary vertices}, the other ones are \emph{mobile vertices}. A hosohedron is called \emph{minimal} if for $ a> 1 $ it does not contain edges that are meridians, i.e. edges connecting the \black--vertex with the \cross--vertex.

Let $ f \colon \Sigma' \rightarrow S^2 $ be an orientation preserving finite-fold branched covering of the sphere by an oriented closed surface $ \Sigma' $. Let $ \Theta (f) $ denote the minimal hosohedron whose vertex set contains all the critical values of the covering $ f $ and any mobile vertex is a critical value (so that stationary vertices may not be critical values). Let $ a_f = a $ be the number of meridians of the hosohedron $ \Theta (f) $. In particular, if $ f $ has no critical values other than poles, then $ a_f = 1 $ and $ \Theta (f) $ consists of two (stationary) vertices and one (arbitrary) meridian.

\emph{The graph} $ \Gamma (f) $  \emph{of branched covering} $ f \colon \Sigma' \rightarrow S^2 $ is an oriented graph on the surface $ \Sigma' $ satisfying the conditions:
\renewcommand{\labelitemi}{$-$}
\begin{itemize}
	\item
 the vertices of $ \Gamma (f) $ are either \emph{stationary} or \emph{mobile};
 \item the stationary vertices are divided into \black-- and \cross--vertices that are the preimages of \black-- and \cross--vertices of  $ \Theta (f) $ respectively;
  \item the mobile vertices are the critical points of   $ f $ other than stationary vertices; 
  \item let $ \Theta ^ 1 \subset S^2 $ be the union of the vertices and edges of the graph $ \Theta (f) $; the edges of  $ \Gamma (f) $ are the components of the preimage $ f ^ {- 1} (\Theta ^ 1) $, cut at the vertices of  $ \Gamma (f) $; each edge is oriented in the same way and marked with the same number as its image in $ \Theta (f) $.
\end{itemize}
In other words, $ \Gamma (f) $ is obtained from $ f ^ {- 1} (\Theta (f)) $ by removing mobile vertices of valency $ 2 $ that are different from the singular points of the surface $ \Sigma' $. (The valency of a vertex is the number of edges adjacent to it.)

It is clear that the valencies of the mobile vertices are even
and not less than $ 4 $,
and the valencies of the stationary ones are multiples of $ a_f $. For example, for an $ n $-fold branched covering $ f \colon S^2 \rightarrow S^2 $ with two critical points lying over the poles, $ a_f = 1 $, as indicated above, and $ \Gamma (f) $ is a hosohedron, not minimal for $ n> 1 $.

\begin{figure} [h]
\begin{center}
\scalebox {0.6} {\includegraphics {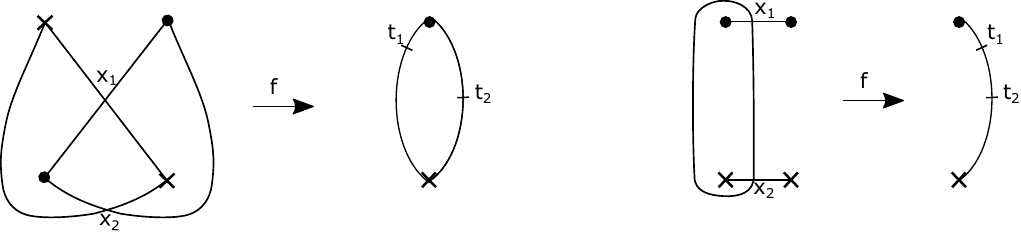}}
\end{center}
\caption {Graphs of two-fold branched coverings.}
\label{MorphDeg2}
\end{figure}

\begin{figure} [h]
\begin{center}
\scalebox {0.6} {\includegraphics {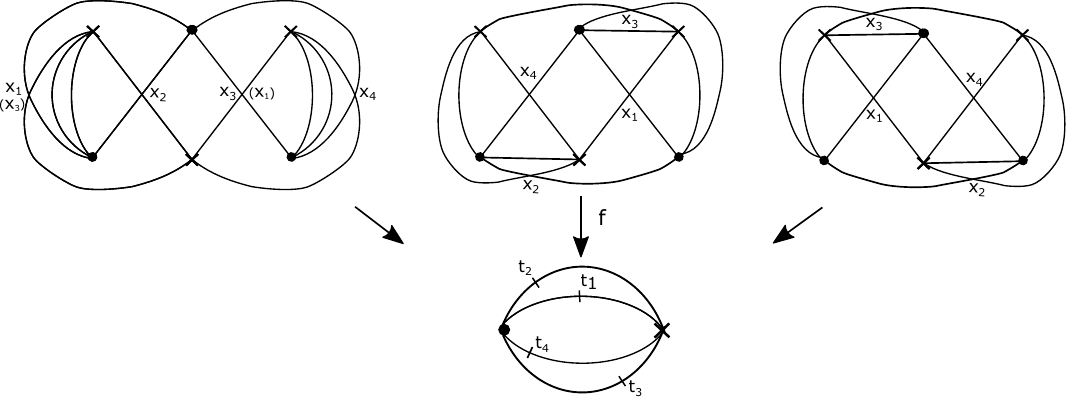}}
\end{center}
\caption {Graphs of three-fold branched coverings.}
\label{MorphDeg3}
\end{figure}

Figures \ref {MorphDeg2}, \ref {MorphDeg3} show examples of graphs $ \Gamma (f) $ and their images $ \Theta (f) $ for branched coverings of degrees $ 2 $ and $ 3 $.
In Figure \ref {MorphDeg3}, the left arrow corresponds to two coverings, turning into each other by permutation of the vertices $ x_1, x_3 $.

It is clear that for isomorphic coverings $ f_1 $, $ f_2 $, the graphs $ \Gamma (f_1) $, $ \Gamma (f_2) $ are isomorphic and $ \Theta (f_1) = \Theta (f_2) $. 
 
\begin{opr}
\label{hoso}
A hosohedral graph with parameter $ a $ \rm {is a finite directed graph $ \Gamma $ embedded into an
	oriented closed surface $ \Sigma' $
	and endowed with the following structure:
	\newcounter {F}
	\begin{list} {(\arabic {F})} {\usecounter {F}}
	
	\item \label{a} the edges of  $ \Gamma $ are
	marked with the
	integers from $ 1 $ to $ a $;
	\item \label{singv} all singular points of $ \Sigma' $ are  vertices of $ \Gamma $ (the \emph{singular} vertices);
	\item the vertices of $ \Gamma $ are divided into \emph{stationary} and
	\emph{mobile};
	each stationary vertex has one of the two \emph{labels}, or \emph{colors}:  \black- \, and \cross-;
	
	\item \label{val} the valency of each stationary vertex is a
	multiple of $ a $, and of each mobile one is
	even;
	
	\item \label{Stat} all edges adjacent to a \black--vertex (\cross--vertex) are outgoing (respectively, incoming) and their labels on each component of  $ \Sigma' $ change cyclically from $ 1 $ to $ a $ when going round this vertex with a positive (respectively, negative) direction, given by the orientation of the surface component;
	
	\item \label{Mob} the labels of all edges adjacent to
	a mobile vertex, are the same and define its \emph{label}, or \emph{color}, and the incoming and outgoing edges alternate on each surface component;
	
	\item \label{Face} each face of  $ \Gamma $ is simply connected and has  one \black-- and  one \cross--vertex on the boundary.
	\end{list}
}
\end{opr}
Note that due to conditions (\ref {singv}) and (\ref {Face}) the graph $ \Gamma $ is connected.

For any face $ F $ of  $ \Gamma $, by conditions (\ref {Face}), (\ref {Stat}) and (\ref {Mob}), \black-- and \cross--vertices split its boundary $ \partial F $  into two oriented paths:
$$
\partial ^ i F, \quad \partial ^ {i + 1} F;
$$
the first path consists of edges with the label $ i $, and the second one with the label $ (i + 1) \mathrm {mod} \, a $.

The union of all (open) edges labeled with $ i $ and the mobile vertices adjacent to them is called the \emph{$ i $-th monochrome part} of  $ \Gamma $. This definition is correct due to condition (\ref {Mob}). A hosohedral graph other than a hosohedron is called \emph{minimal} if it does not contain mobile vertices of valency $ 2 $ and for any $ i = 1, \ldots, a $ its $ i $-th monochrome part contains a mobile vertex. A path in  $ \Gamma $ is called \emph{monochrome} if it lies in a monochrome part, i.e. it does not go through stationary vertices and all its edges are marked with the same label. Let us define a binary relation $ \prec $ on the set of mobile vertices of the graph: $ u \prec v $, if there is a monochrome oriented path going from $ u $ to $ v $. A hosohedral graph is called \emph{admissible} if the relation $ \prec $ is a partial order. It's clear that
the relation $ \prec $ is transitive, so it is a partial order if and only if the graph $ \Gamma $ does not contain oriented monochrome cycles.

\begin{figure} [h]
\begin{center}
\scalebox {0.6} {\includegraphics {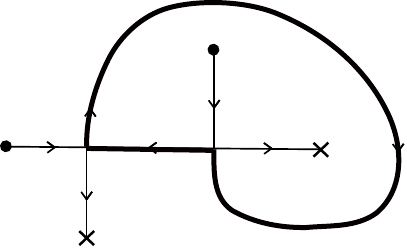}}
\end{center}
\caption {Inadmissible hosohedral graph.}
\label{nonadm1}
\end{figure}

The simplest example of an admissible hosohedral graph is a hosohedron. An example of an inadmissible  hosohedral graph is shown in Figure \ref {nonadm1} (the bold line indicates an oriented monochrome cycle).
A minimal admissible hosohedral graph (with parameter $ a $) is called a \emph{hosograph} (respectively, $ a$-\emph{hosograph}).

Let $ \Gamma_1 $, $ \Gamma_2 $ be hosohedral graphs. The mapping $ \varphi \colon \Gamma_1 \rightarrow \Gamma_2 $ is called a \emph{morphism} of hosohedral graphs if it preserves the hosohedral structure, i.e.  takes the vertices, edges, faces of $ \Gamma_1 $, respectively, to vertices, edges, faces of $ \Gamma_2 $, keeping the type of vertices, orientations and edge labels.	
\begin{te}
\label{graphf}
\begin{enumerate}
\item A hosohedral graph $ \Gamma $ is admissible if and only if  there is a morphism of this graph onto a hosohedron.
\item
A graph~$ \Gamma $ can be represented by the graph
of branched covering of the sphere if and only if it is a hosograph.
\item \label{3}
For a given hosograph $ \Gamma $, the branched covering $ f $ with $ \Gamma (f) = \Gamma $ is defined
up to isomorphism by the images of the mobile vertices of $ \Gamma, $ or, equivalently,  by the critical values of $ f $ with the sets of their preimages.
\end{enumerate}
\end{te}
\begin{proof}
To prove the first statement, we use the proof idea of a similar theorem \cite [Theorem 4.11] {Degt}.

Let $ \Gamma $ be a hosohedral graph with parameter $ a $.

We construct a morphism $ f \colon \Gamma \rightarrow \Theta $, where $ \Theta $ is a hosohedron with $ a $ meridians, the mobile vertices of $ \Gamma $ being specified during construction. The morphism takes stationary vertices to stationary ones. The distance from the \black--vertex of the hosohedron to a point on its meridian defines the linear order < on each meridian. Let us extend the partial ordering $ \prec $ on the vertex set of the $ i $-th monochrome part of $ \Gamma $ to an arbitrary linear order and map this set monotonically onto an increasing set of equally spaced points of the $ i $-th meridian of  $ \Theta $ , thus specifying the mobile vertices of the hosohedron. It is clear that the resulting mapping is extended to the entire $ i $-th monochrome part.
Thus, for any face $ F $ of  $ \Gamma $, the mapping $ f $ is defined on both parts $ \partial ^ i F $, $ \partial ^ {i + 1} F $ of its boundary. Since $ F $ is simply connected, $ f $ is extended to this face, taking it onto a 2-gon between the $ i $-th and $ (i + 1) $-th meridians on the sphere. This completes the construction of  $ f $.

Conversely, the linear order < on the $ i $-th meridian of $ \Theta $ induces, via the morphism $ f \colon \Gamma \rightarrow \Theta $, a partial order on the vertex set of the $ i $-th monochrome part of  $ \Gamma $. Therefore, $ \Gamma $ is admissible.

Let us prove the second statement.

Let $ f $ be a branched covering of the sphere. Since the faces of  $ \Theta (f) $ do not contain the critical values of $ f $, the faces of  $ \Gamma (f) $ are simply connected. The other conditions of Definition \ref{hoso} are satisfied for $ \Gamma (f) $ as it  immediately follows from the definition of $ \Gamma (f) $. Since the mapping $ f \colon \Gamma (f) \rightarrow \Theta (f) $ preserves the orientation of the edges, the relation $ \prec $ on the $ i $-th monochrome part of $ \Gamma (f) $ is induced by the linear order on the $ i$-th meridian of $ \Theta (f) $ and therefore is a partial order. Thus,  $ \Gamma (f) $ is admissible and $ f $ is a morphism. Moreover, it is clear that $ \Gamma (f) $ is minimal and therefore is a hosograph.

Conversely, let $ \Gamma $ be a hosograph on a surface $ \Sigma' $ and $ f \colon \Gamma \rightarrow \Theta $ be the morphism constructed in the proof of the first statement. We may assume that~$ f $ is induced by some continuous mapping $ \Sigma' \rightarrow S^2 $, which we also denote by $ f $.
It is clear that the points of the faces and the interior points of the edges are regular points of the mapping $ f $. Therefore, all its critical points are isolated, and $ f $ is a branched covering.

Let us prove the third statement. From the proofs of the first and second statements it follows that there is a covering $ f $ with $ \Gamma (f) = \Gamma $ and with given images of the mobile vertices of  $ \Gamma $. An enumeration $ t_1, \ldots, t_k $ of
the sphere poles and these images defines a \emph{constellation}, i.e. a sequence of permutations generating the monodromy group of the covering $ f $ (see \cite [1.2.3] {LZ}). According to \cite [\S 1.6] {LZ}, the graph~$ \Gamma $ together with the points $ t_1, \ldots, t_k $ uniquely determine this constellation. Finally, \cite [Proposition 1.2.16] {LZ} makes the branched coverings of the sphere with the given constellation isomorphic.
\end{proof}

For branched covering $$
f \colon \Sigma' \rightarrow S^2
$$
the same letter  denote the corresponding morphism $$
f \colon \Gamma (f) \rightarrow \Theta (f).
$$

Note that there are coverings $ f_1, \! f_2 $ with $ \Gamma (f_1) \! \cong \! \Gamma (f_2) $ \! and $ {\Theta (f_1) \! \ncong \! \Theta (f_2)} $, since the images of two non-comparable mobile vertices of the same color can either coincide or be  different. The smallest degree of such coverings is $ 4 $.

\begin{za}
\label{basegraph}
In the choice of the base graph on the sphere used in the definition of graph of branched covering, other options are possible.
One can construct a "latitudinal" version of the theory, replacing the hosohedron with a base graph, whose vertices and edges lie on one meridian and several parallels. Such a graph, like a hosohedron, is associated with polar coordinates on the plane that is the image of the sphere under the stereographic projection from a pole. The rectangular coordinates on this plane are associated with a base graph formed by a wedge of circles on the sphere touching each other at the pole. Using such a base graph
one can, for example, describe the cellular decomposition of the space $ N (\Sigma, n) $ built in \cite {DE}.
\end{za}

\subsection {Graph transformations}
\label{Transform}

For $ i = 0,1 $ let $ \Gamma_i $ be a hosograph on a surface $ \Sigma_i $ and $ f_i $ be a branched covering with $ \Gamma (f_i) = \Gamma_i $.
The hosograph~$ \Gamma_1 $ is a \emph{perturbation} of~$ \Gamma_0 $, and  $ \Gamma_0 $ is a \emph{degeneration} of $ \Gamma_1 $, if $ f_1 $ is a perturbation of  $ f_0 $.

It's clear that
if $ \Gamma_i $ is $ a_i $-hosograph, then $ a_1 \geq a_0 $.

The operations of transition from $ \Gamma_0 $ to $ \Gamma_1 $ and back are also called a \emph{perturbation} and a \emph{degeneration of hosograph}.

The next three degenerations
are called \emph{elementary}:
\newcounter {E}
\begin{list} {(\alph {E})} {\usecounter {E}}
\item $ mob (v_1, v_2) $ degeneration that contracts to a mobile vertex $ v $ all parallel edges connecting two different adjacent mobile vertices $ v_1, v_2 $;
moreover, $ v $ is non-singular if and only if $ v_1, v_2 $ are connected by a single edge;
\item $ stat (v) $ degeneration such that, for all edges connecting one or more stationary vertices of the same color with a mobile vertex $ v $, it contracts them to a stationary vertex; moreover, the resulting stationary vertex is non-singular if and only if $ v $ is connected to each of the indicated stationary vertices by a single edge;

\item $ edg $ degeneration that cuts each face $ F $ between the $ i $-th and $ (i + 1) $-th monochrome parts of  $ \Gamma_1 $ and glues $ \partial ^ i F $ to $ \partial ^ {i + 1} F $ so that no two vertices are glued into one.

\end{list}

The perturbations inverse to elementary degenerations are called \emph{elementary perturbations.}

It is clear that elementary transformations take hosographs to hosographs, therefore we can talk about elementary transformations of coverings and, moreover, their isomorphism classes.

\begin{pr}
\label{stat}
Let $ f_1 $ be a perturbation of a covering $ f_0, $ having a perturbation domain $ V $ and obtained using $ stat ^ {- 1} (v_1), $ where $ v_1 $ is a critical point of  index~$ 2 $ \textup (i.e. $ f_1 (v_1) $ is a simple critical value \textup). Then $ f_1 $ is unique up to $ (V, w) $-equivalence \emph{(see Definition in Subsection \ref {Deform})}, where $ w \in V $ is the corresponding stationary vertex of  $ \Gamma (f_0) $.
\end{pr}
\begin{proof}
Let $f_2 \colon \Sigma_2 \rightarrow S^2 $ be a perturbation of  $f_0 $ that has the same perturbation domain $ V $ and is obtained using $ stat ^ {- 1} (v_2) $, where $ v_2 $ is a critical point of  index~$ 2 $. Let us prove that $ f_1, f_2 $ are $ (V, w) $\nobreakdash-equivalent.

In the (closed) disk $ V = \bar {U} _w $ (see the notation in Subsection \ref {Deform}), choose a closed sector $ D $ containing the points $ f_1 (v_1) $ and $ f_2 (v_2) $. Since $ v_1 $, $ v_2 $ are critical points of index $ 2 $, on both surfaces $ \Sigma_1, \Sigma_2 $ the surface~$ \widetilde {U} _v $ is a disk if the point $ w $ is non-singular, and a ring otherwise. Therefore, we can assume that $ \Sigma_1= \Sigma_2 $ is  the same surface. Outside $ \widetilde {U} _v $, the coverings $ f_1, f_2 $ coincide. In $ \widetilde {U} _v $ on the preimages of the boundary radii of the sector $ D $, these coverings are obviously homotopic. This homotopy can be extended to the preimages of the sector $ D $ and its complements in $ V $, since on all components of these preimages, except one, the coverings $ f_1, f_2 $ are univalent, and on the remaining component they are   double-coverings of $ D $ with a single branch point.\qed
\end{proof}

\section {Cell structure}
\label{Cell}
\subsection {$ 0 $-in-graph}
We   consider the hosohedron $ \Theta $ as a subset of the sphere $ S^2 $ so that we can talk about the angles between the meridians and the distances from the poles to the vertices of the hosohedron. Select a meridian on the sphere, which we call the \emph{prime meridian}. The union of $ \Theta $ with the prime meridian is called the $ 0 $-in-\emph{hosohedron}. In particular, we assume that for a minimal $ 0 $-in-hosohedron without mobile vertices, the prime meridian is its only edge. Let $ \Gamma $ be a hosohedral graph with parameter $ a + 1 $ whose edges are labeled with the integers from $ 0 $ to $ a $. A graph $ \Gamma $ is called a $ 0 $-in-\emph{graph} (more precisely, $ a$-$ 0 $-in-\emph{graph}) if there exists a morphism $ f \colon \Gamma \rightarrow \Theta $ that takes the $ 0 $-th monochrome part to the prime meridian, $ \Gamma $ does not contain mobile vertices of valency $ 2 $ and for any $ i = 1, \ldots, a $ its $ i $-th monochrome part contains a mobile vertex. Thus, a $ 0 $-in-graph is a hosograph if and only if its $ 0 $-th monochrome part contains at least one (mobile) vertex; if its $ 0 $-th monochrome part does not contain vertices, it will be a hosograph after the part is removed.

Let $ \Gamma $ be a $ 0 $-in-graph. We call it $ 0 $-in-\emph{graph of branched covering}~$ f $ if $ \Gamma $ is the union of $ \Gamma(f) $ with the preimage of the prime meridian. Perturbations/ degenerations of $ 0 $-in-graphs are defined using coverings in the same way as the perturbations/\-degenerations of hosographs.
\subsection {Mapping L}
For an $ a$-$ 0 $-in-graph $ \Gamma $ on a surface $ \Sigma' $, we denote by $ E _ {\Gamma} $ the set of isomorphism classes of branched coverings 
$ f \colon \Sigma' \rightarrow S^2 $ such that the $ 0 $-in-graph of $ f $ is isomorphic to $ \Gamma $.
For each $ i = 0,1, \ldots, a $ we enumerate the mobile vertices of the $ i $-th monochrome part of $ \Gamma $, denoting them by $ v_ {ij} $. Let $ \nu_i = \# \{v_ {ij} \} $ be the number of such vertices and $ k = \sum_i \nu_i $.
	Let $ L _ {\Gamma} \colon E _ {\Gamma} \rightarrow (S^2) ^ k $ be the mapping that takes the covering $ f $ to the set $ f (v_ {ij}) $ of its critical values other than poles.
	
	If $ \Gamma $ has no mobile vertices, then $ a = 0 $ and the mapping $ L _ {\Gamma} $ is undefined.
	
	We denote by $ \sigma \colon [0,2 \pi] \times [0, \pi] \rightarrow S^2 $ a parametrization of the sphere such that $ \sigma (\lambda, \varphi) $ is a point with the longitude $ \lambda $ and the latitude $ \varphi $,  $ \sigma (0 \times [0, \pi]) = \sigma (2 \pi \times [0, \pi]) $ is the prime meridian, and $ \sigma (0,0) $,  $ \sigma (0, \pi) $ are the \black-- and \cross--vertices, respectively. Let $ \sigma ^ k \colon ([0,2 \pi] \times [0, \pi]) ^ k \rightarrow (S^2) ^ k $ be the power of the parameterization. It is clear that the restriction of the mapping $ \sigma ^ k $ to $ ((0,2 \pi) \times (0, \pi)) ^ k $ is a topological embedding.
	We denote by $ \lambda_ {ij} \in [0,2 \pi], \varphi_ {ij} \in [0, \pi] $ the coordinates in $ ([0,2 \pi] \times [0, \pi ]) ^ k $, where $ i = 0,1, \ldots, a, j = 1, \ldots, \nu_i $.
	\begin{te}
	\label{mapL}
	The mapping $ L _ {\Gamma} $ is a topological embedding. The 
	intersection $ \mathring {M} = (\sigma ^ k) ^ {- 1} (L _ {\Gamma} (E _ {\Gamma})) \cap ([0,2 \pi) \times [0, \pi ]) ^ k $ is an open convex polytope of dimension $ k + a $.
	\end{te}
	\begin{proof}
	If $ L _ {\Gamma} ([f_1]) = L _ {\Gamma} ([f_2]) $, then the sets of critical values of the coverings $ f_1 $, $ f_2 $ coincide and the graphs $ \Gamma (f_1) $, $ \Gamma (f_2) $ are isomorphic. Therefore, $ [f_1] = [f_2] $ due to  (\ref{3}) of Theorem \ref {graphf}. So $ L _ {\Gamma} $ is injective.
	
	For $ [f] \in E _ {\Gamma} $, the interior of the perturbation domain $ V $ of the covering $ f $ obviously determines a neighborhood $ \widetilde {V} $ of the point $ x = L _ {\Gamma} ([f]) $. According to
	Theorem \ref {graphf} the class~$ [f] $ is uniquely determined by the graph $ \Gamma $ and the set of points $ f (v_ {ij}) $. Therefore
	$$
	L _ {\Gamma} (E _ {\Gamma} \cap U _ {[f], V}) = \widetilde {V} \cap L _ {\Gamma} (E _ {\Gamma}).
	$$
	Conversely, any neighborhood $ \widetilde {V} $ of the point $ x $ gives a perturbation domain $ V $ of  $ f $ with $ L _ {\Gamma} (E _ {\Gamma} \cap U _ {[f], V}) = \widetilde {V} \cap L _ {\Gamma} (E _ {\Gamma}) $. Taking into account Proposition \ref {base}, we obtain that $ L _ {\Gamma} $ establishes a bijection between the bases of  topologies at the point~$ [f] $ and at the point $ x $. By virtue of the above, $ L _ {\Gamma} $ is injective, therefore it is a topological embedding.
	
	From the proof of Theorem \ref {graphf} it follows that each class $ [f] \in E _ {\Gamma} $ is uniquely determined by the following data: 1) by   $ a $   meridians of the hosohedron $ \Theta (f) $ provided that if the $ i $-th meridian is obtained from the prime one
	by turning through an angle $ x_i $, then $ 0 <x_1 <x_2 <\ldots <x_a <2 \pi $, and 2) by the images of the mobile vertices of $ \Gamma $ satisfying the condition $ v_ {ij} \prec v_ {ik} \Leftrightarrow x_ {ij} <x_ {ik} $, where $ x_ {ij} $ is the latitude of  $ f (v_ {ij}) $.
	Therefore, taking into account that the restriction of  $ \sigma ^ k $ to $ ((0,2 \pi) \times (0, \pi)) ^ k $, as noted above,
	is a topological embedding, we see that $ \mathring {M} $ is defined in $ ([0,2 \pi] \times [0, \pi]) ^ k $ by the equations $ \lambda_ {01} = \lambda_ {02} = \ldots = \lambda_ {0 \nu_0} = 0 $, $ \lambda_ {i1} = \lambda_ {i2} = \ldots = \lambda_ {i \nu_i}, i = 1, \ldots, a $, and by specified above  inequalities for the longitudes and latitudes of the images of the mobile vertices of  $ \Gamma $. Therefore, $ \mathring {M} $ is an open convex polytope of dimension $ k + a $.
	\end{proof}
	\begin{sle}
	\label{OpenCell}
	The space $ \bar {X} _ {\Sigma, n} $ is divided into open cells $ E _ {\Gamma} $. For an $ a$-$ 0 $-in-graph $ \Gamma $, the dimension of a cell $ E _ {\Gamma} $ is $ \nu_0 + \sum_ {i = 1} ^ a (1+ \nu_i) = k + a, $ where $ \nu_i $ is the number of mobile vertices of the $ i $-th monochrome part of  $ \Gamma $. \qed
	\end{sle}
	\begin{za}
	\label{coord}
	
	From the proof of Theorem \ref {mapL} it follows that the polytope $ \mathring {M} \subset ([0,2 \pi] \times [0, \pi]) ^ k $ is determined by the equations and  inequalities, indicated in the theorem. 
	 Therefore, as coordinates in $ \mathring {M} $ one can take the longitudes $ x_1 = \lambda_ {11}, x_2 = \lambda_ {21}, \ldots, x_a = \lambda_ {a1} $ and the latitudes $ x_ {ij} = \varphi_ {ij} $ of the images of the mobile vertices $ v_ {ij} $ of $ \Gamma $, where $ i = 0,1, \ldots, a $, $ j = 1, \ldots, \nu_i $.
	
	\end{za}
\subsection {Characteristic mapping of a cell}
By Theorem \ref {mapL}, the set $ \mathring {M} $ is an open convex polytope and, therefore, is homeomorphic to the open ball $ \Int D ^ {\dim E _ {\Gamma}} $.
We denote by $ M $ the closure of this set in $ ([0,2 \pi] \times [0, \pi]) ^ k $. It is clear that $ M \cong D ^ {\dim E _ {\Gamma}} $ is a bounded closed polytope.
\begin{te}
\label{Clcell}
Let $ \bar {E} _ {\Gamma} $ denote the union of $ E _ {\Gamma} $ with the cells corresponding to all degenerations of the $ 0 $-in-graph $ \Gamma $. Then
\begin{enumerate}
\item\label{3.1} $ \bar {E} _ {\Gamma} $ is the closure of the cell $ E _ {\Gamma}; $
\item there is a characteristic mapping $ \chi _ {\Gamma} \colon M \rightarrow \bar {X} _ {\Sigma, n} $ of $ E _ {\Gamma}, $ i.e. $ \chi_ { \Gamma} $ takes $ \mathring {M} $ homeomorphically to $ E _ {\Gamma}, $ and $ \partial M $ to the union of cells of dimension $ <\dim E _ {\Gamma}; $
\item\label{3.3} $ \chi _ {\Gamma} (M) = \bar {E} _ {\Gamma} $.
\end{enumerate}
\end{te}
\begin{proof}
\begin{enumerate}
\item Let $ [f_0] $ be a degeneration of a class $ [f_1] \in E _ {\Gamma} $ and $ V $ be the corresponding perturbation domain of $ f_0 $. For any other perturbation domain $ V_1 $ of $ f_0 $, the composition $ \psi \circ f_1 $, where $ \psi $ is a homeomorphism of the sphere onto itself, taking $ V $ to $ V_1 $, is a perturbation of  $ f_0 $ with  $[\psi \circ f_1]\in E _ {\Gamma} $  corresponding to the perturbation domain  $ V_1 $. Since the collection $ \{U _ {[f_0], V} \} _ V $ forms the base of topology at the point $ [f_0] $, this point lies in the closure of the cell $ E _ {\Gamma} $.
	
	Conversely, let $ [f_0] $ be an adherent 
	 point of the set $ E _ {\Gamma} $. Then for some perturbation domain $ V $ of $ f_0 $ the neighborhood $ U _ {[f_0], V} $ of this point intersects $ E _ {\Gamma} $. Therefore, $ [f_0] \in \bar {E} _ {\Gamma} $.
	
	\item
	Put $ \chi _ {\Gamma} | _ {\mathring {M}} = L _ {\Gamma} ^ {- 1} \circ \sigma ^ k \colon \mathring {M} \rightarrow \bar {X} _ { \Sigma, n} $ and extend $ \chi _ {\Gamma} $ to $ \partial M $ by inverse induction on the dimension of the faces of the polytope, proving that the image of each open face is a cell.
	Let  $ \chi _ {\Gamma} $ be extended to an open face $ F_1 $ and $ \chi _ {\Gamma} (F_1) = E _ {\Gamma_1} $ is a cell. Take an open face $ F_0 $ lying on the boundary of $ F_1 $ with $ \dim F_0 = \dim F_1-1 $, and a point $ x \in F_0 $. To define the class $ [f_0] = \chi _ {\Gamma} (x) $,
	it suffices to indicate, according to Theorem \ref {graphf}, the
	graph $ \Gamma_0 = \Gamma (f_0) $ and the images of its mobile vertices under the mapping~$ f_0 $.
	The following cases are possible.
	\begin{enumerate}
	\item The face $ F_0 $ lies on the boundary of the half-space $ x_ {ij} \leq x_ {ik} $ (see the notation of coordinates in Remark \ref {coord}). It follows from the description of the image of the mapping $ L _ {\Gamma} $ (see the proof of Theorem \ref {mapL}) that there exist a class $ [f_1] \in E _ {\Gamma_1} $ and mobile vertices $ v_ {ij}, v_ { ik} $ of  $ \Gamma_1 $ such that $ f_1 (v_ {ij}) $, $ f_1 (v_ {ik}) $ are neighboring vertices of the hosohedron $ \Theta (f_1) $, and the latitudes of the latter satisfy the inequality $ x_ { ij} <x_ {ik} $. Therefore, $ v_ {ij}, v_ {ik} $ are adjacent and $ v_ {ij} \prec v_ {ik} $.
	In this case, the graph $ \Gamma_0 $ is obtained from  $ \Gamma_1 $ by the elementary degeneration $ mob (v_ {ij}, v_ {ik}) $, contracting all parallel edges that connect $ v_ {ij} $ with $ v_ {ik} $.
	\item The face $ F_0 $ is defined by the equation $ x_ {i1} = 0 $ or $ x_ {i \nu_i} = \pi $. Then there is a class $ [f_1] \in E _ {\Gamma_1} $ and a mobile vertex $ v $ of  $ \Gamma_1 $ such that the vertex $ f_1 (v) $ is adjacent to a pole in $ \Theta (f_1) $. Therefore, $ v $ is connected by edges either with one or with several stationary vertices of the same color. Thus,  $ \Gamma_0 $ is obtained from  $ \Gamma_1 $ by the elementary degeneration $ stat (v) $, which contracts all edges connecting these stationary vertices with $ v $ to a stationary vertex. Moreover, if $ \nu_i = 1 $, then the $ i $-th monochrome part in the resulting graph is removed.
	
	\item The face $ F_0 $ is defined by the equation $ x_ {1} = 0 $, $ x_ {i} = x_ {i + 1} $ or $ x_a = 2 \pi $. Then $ \Gamma_0 $ is obtained from  $ \Gamma_1 $ by an elementary degeneration $ edg $.
	
	In all these cases, the elementary degenerations and the coordinates of the point $ x $ uniquely determine the images of the mobile vertices of  $ \Gamma_0 $.
	
	\end{enumerate}
	Let us prove that $ \chi _ {\Gamma} $ is continuous at any point $ x \in \partial M $ (and then the constructed extension does not depend on the choice of the class $ [f_1] \in E _ {\Gamma_1} $ in the above cases). Let $ \chi _ {\Gamma} (x) = [f_0] $. The interior of a perturbation domain $ V $ of $ f_0 $ obviously defines a neighborhood $ \widetilde {V} $ of $ L _ {\Gamma} ([f_0]) $.
	It is clear that $ \chi _ {\Gamma} ^ {- 1} (\bar {E} _ {\Gamma} \cap U _ {[f_0], V}) = M \cap (\sigma ^ k) ^ {- 1} (\widetilde {V}) $. Since the mapping $ \sigma ^ k $ is continuous, the set $ M \cap (\sigma ^ k) ^ {- 1} (\widetilde {V}) $ is open in $ M $. It remains to note that $ \{\bar {E} _ {\Gamma} \cap U _ {[f_0], V} \} _ V $ is a neighborhood base at  $ [f_0] $ according to Proposition \ref {base}. Therefore, $ \chi _ {\Gamma} $ is continuous at  $ x $.
		\end{enumerate}
		Statement (iii) immediately follows from (i) and (ii).
		\end{proof}
A face $ F $ of the polytope $ M $ from the proof of Theorem \ref{Clcell} is called \emph{exceptional} if it is given by the equations $ x_ {i1} = x_ {i2} = \ldots = x_ {i \nu_i} = 0 $ or $ x_ {i1 } = x_ {i2} = \ldots = x_ {i \nu_i} = \pi $, where $ i \in A \subset \{1,2, \ldots, a \}, A \neq \varnothing $. Eliminate from $ ([0,2 \pi] \times [0, \pi]) ^ k $ the factors corresponding to the coordinates
	$$
	\lambda_ {ij}, \quad \varphi_ {ij}, \quad i \in A, j = 1, \ldots, \nu_i $$
	and denote by $ \prr_A $ the projection of the product $ ([0,2 \pi] \times [0, \pi]) ^ k $ onto the resulting space. It is clear that each exceptional face corresponds to an admissible hosohedral graph in which the $ i $-th monochrome parts for $ i \in A $ do not contain mobile vertices, and after removing these parts, we get a $ 0 $-in-graph. Elementary degenerations of such graphs are defined in the same way as elementary degenerations of hosographs.
	\begin{pr}
	\label{prF}
	For every exceptional face $ F $ of the polytope $ M $ there are non-exceptional  faces of $ F $ whose images under the projection $ \prr_A $ are equal to $ \prr_A (F) $. These are those and only those faces \textup whose $ 0 $-in-graph is obtained from the above-mentioned hosohedral graph of the face $ F $ by removing the $ i $-th monochrome parts for $ i \in A $.
	\end{pr}
	\begin{proof}
	For the proof, it suffices to apply the elementary degenerations $ edg $ to all $ i $-th monochrome parts of $ \Gamma $, $ i \in A $,  corresponding to the face $ F $.
	\end{proof}
	\begin{za}
	\label{restr}
	Let $ F $ be any face of the polytope $ M $, $ \chi $ be the restriction of the mapping $ \chi _ {\Gamma} $ to $ F $ if $ F $ is non-exceptional, and $\chi= \chi _ {\Gamma} \circ \prr_A ^ {- 1} \colon \prr_A (F) \rightarrow \bar {X} _ {\Sigma, n} $ if $ F $ is  exceptional. It obviously follows from the proof of Theorem \ref{Clcell} that $ \chi $ is a characteristic mapping of some cell lying on the boundary of $ E _ {\Gamma} $. Moreover, for an exceptional face $ F $, we have $ \dim \prr_A (F) = \dim F- | A | $, where $ | A | $ is the number of elements in the set $ A $.
	\end{za}
	The following statement obviously follows from the proof of Theorem \ref{Clcell}.
	\begin{sle}
	\label{el}
	Any perturbation/\-degeneration of a hosograph is a composition of elementary perturbations/\-degenerations. \qed
	\end{sle}
	\begin{sle}
	\label{compact}
	The space $ \bar {X} _ {\Sigma, n} $ is a finite cellular space and therefore compact.
	\end{sle}
	\begin{proof}
	It is clear that for a fixed $ n $ the number of graphs of branched coverings of degree $ n $ is finite.
	 \end{proof}
	Corollary \ref {compact} gives a  proof that  the space $ \bar {X} _ {\Sigma, n} $ is compact, independent of \cite [Theorem 2.5] {NT}.
	
	We fix distinct points $ w_1, \ldots, w_r \in S^2 $ and partitions $ \lambda_1, \ldots, \lambda_r $ of the number $ n $. Put $ w = (w_1, \ldots, w_r), \lambda = (\lambda_1, \ldots, \lambda_r) $. Let $ H (\Sigma, n, w, \lambda) \subset X _ {\Sigma, n} $ be the set of classes of branched coverings $ f $ such that the critical values of $ f $ lying in $ S^2 \setminus \{w_1 , \ldots, w_r \} $ are simple, and $ \lambda_i $ is the collection of local multiplicities of  $ f $ at the points of $ f ^ {- 1} (w_i) $. For $ r = 0 $, we have $ H (\Sigma, n, w, \lambda) = H (\Sigma, n) $. In general case the space $ H (\Sigma, n, w, \lambda) $ is  not connected (see, in particular, \cite [Example 5.5.8] {LZ}). All coverings in $ H (\Sigma, n, w, \lambda) $ have the same collections of multiplicities $ k_1, \ldots, k_r $ of the critical values $ w_1, \ldots, w_r $.
		
		Let $ H_w $ be any connected component of the space $ H (\Sigma, n, w, \lambda) $.
		Due to (\ref{3.1}), (\ref{3.3})  of Theorem \ref {Clcell} we obtain the following statement.
		\begin{sle}
		\label{cellsubspace}
		The closure $ \bar {H} _w $ of the set $ H_w $ in $ \bar {X} _ {\Sigma, n} $
		is a cellular subspace, and $ \bar {H} _w \smallsetminus H_w $ is a cellular subspace of $ \bar {H} _w $. \qed
		\end{sle}
		For $ [f] \in \bar {H} _w $ we include the critical values $ w_1, \ldots, w_r $ in the set of stationary vertices of $ \Theta(f) $, and their preimages in the set of stationary vertices of  $ \Gamma (f) $. Let us call the prime meridian and the meridians passing through $ w_1, \ldots, w_r $, \emph{stationary}. 
		Note that as in Section \ref {Transform}, the transformation $ stat (v) $ 
		is also applicable to the new stationary vertices.

\subsection {The cellular structure of the space $ \bar {X} _ {S^2,2} $} \label{ex2}
\begin{figure} [ht]
\begin{center}
\scalebox {0.6} {\includegraphics {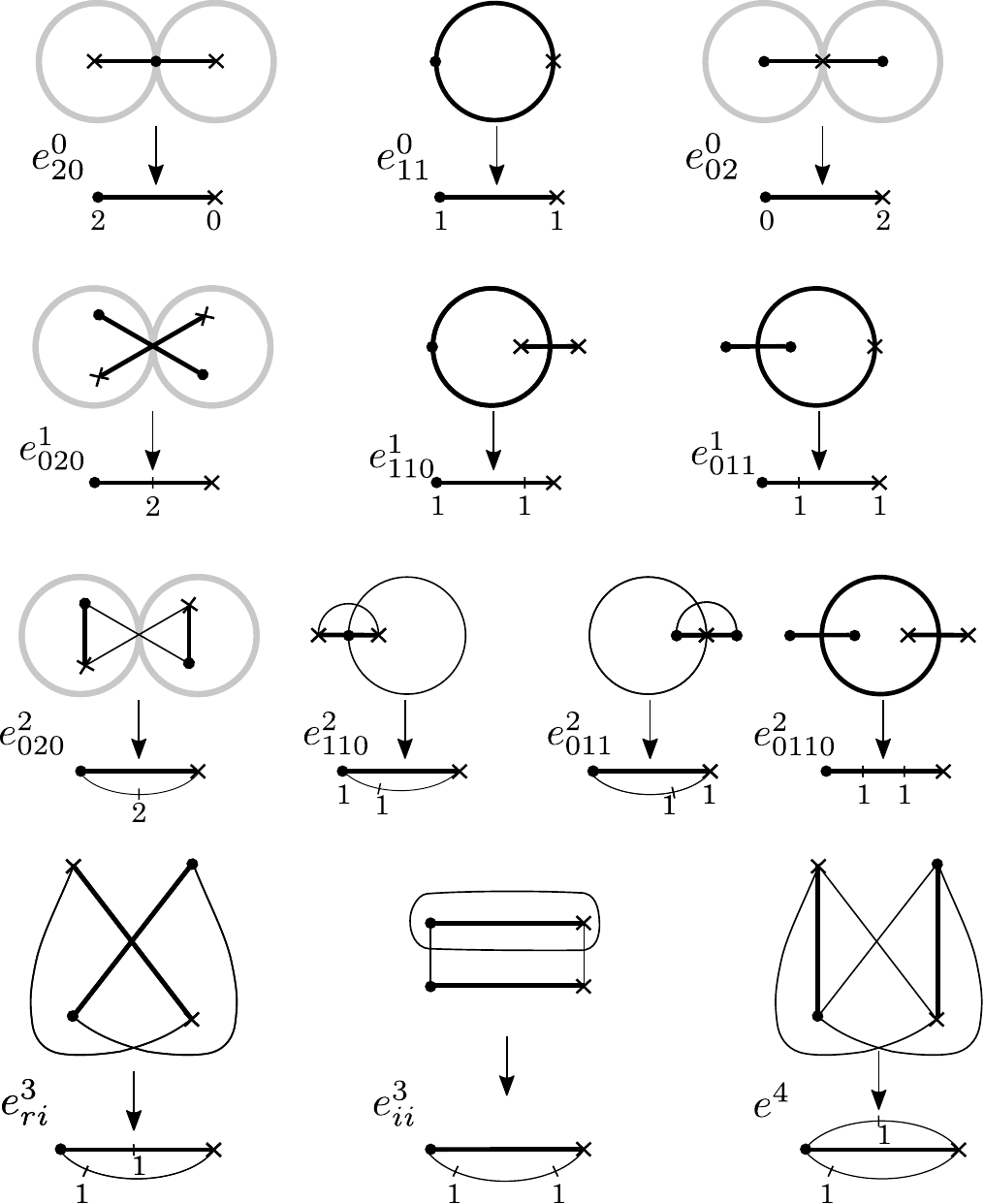}} \\
\end{center}
\caption {Cells in space $ \bar {X} _ {S^2,2} $.} \label{cells}
\end{figure}

Consider the simplest case of two-fold branched coverings of the sphere by the sphere and, as an example, specify  a cellular decomposition of the space $ \bar {X} _ {S^2,2} $. We denote by $ P ^ 2 \cong \mathbb {C} P ^ 2 $ the symmetric square of the sphere $ S^2 $, i.e.  the quotient space of $ (S^2) ^ 2 $ by the permutation of the coordinates. Let  $ LL \colon \bar {X} _ {S^2,2} \rightarrow P ^ 2 $ maps the isomorphism class $ [f] $ of a covering $ f $ to the set of the critical values of $ f $, taken with their multiplicities. According to \cite [Subsection 3.7] {NT} the mapping $ LL $ is a homeomorphism. Therefore, the cellular structure of  $ \bar {X} _ {S^2,2} $ can be described in the language of hosohedra. But for the sake of completeness, we also indicate the corresponding hosographs (see Figure \ref{cells}). In the figure, the cells are denoted by $ e ^ i_ {cr} $, where $ i $ is the dimension of the cell, and $ cr $ characterizes the set of multiplicities of critical values; gray lines represent wedges of spheres, bold lines the prime meridian and its preimage, thin lines  the other meridians and their preimages. The numbers below meridians indicate the multiplicities of the critical values.

The two-dimensional skeleton consists of three spheres touching each other, which are  the closures of the cells $ e ^ 2_ {020}, e ^ 2_ {110}, e ^ 2_ {011} $, and of the plane curvilinear triangle $ e ^ 2_ {0110} $. The touch points are $ 0 $-cells; they are connected by $ 1 $-cells that are the shortest arcs on the spheres and the sides of the triangle $ e ^ 2_ {0110} $.

The cell $ e ^ 3_ {ri} $ is adjacent to $ e ^ 2_ {110}, e ^ 2_ {011} $ and, from four sides, to $ e ^ 2_ {0110} $; the cell $ e ^ 3_ {ii} $ is adjacent to $ e ^ 2_ {020}, e ^ 2_ {110}, e ^ 2_ {011} $ and, from two sides, to $ e ^ 2_ {0110} $. The cell $ e ^ 4 $ is adjacent from two sides to both $ e ^ 3_ {ri} $ and $ e ^ 3_ {ii} $.

Subspace $ \bar {X} _ {S^2,2} \setminus X_ {S^2,2} = e ^ 2_ {020} \cup e ^ 1_ {020} \cup e ^ 0_ {20} \cup e ^ 0_ {02} $ is the closure of the cell $ e ^ 2_ {020} $ and consists of the isomorphism classes of coverings whose definition domain  is a singular surface  (a wedge of two spheres).

For comparison, a more economical cellular decomposition of  $ \bar {X} _ {S^2,2} $ constructed in \cite {DE} consists of a $ 0 $-cell, two $ 2 $-cells, $ 3 $-cells and $ 4 $-cells. The graphs of the corresponding coverings lie on the sphere or on the wedge of two spheres and are the preimages of $ a $ vertical lines on $ \mathbb {C} \cup \infty $, $ a = 0,1,2 $, with the lying on them critical values of covering.
\subsection {Digression: a weakly polyhedral cellular space} \label{def}
The notation $ F <M $  means that the polytope $ F $ is a proper face of the polytope $ M $.
\begin{opr}
\label{wp}

Let $ X $ be a framed cellular space, i.e.  a cellular space together with a characteristic mapping $ \chi_E \colon D_E \rightarrow X $ chosen for each (closed) cell $ E $. We   assume that
\begin{enumerate}
\item \label{wp1}
$ D_E $ is a compact convex polytope of dimension $ d = \dim E $;
\item \label{wp2} for any
face $ F <D_E $, there exist a cell $ E_1 \subset \partial E $ of dimension $ \leq \dim F $ (maybe the same for different proper faces of $ D_E $) and a linear surjective mapping $ p_ {E, F} \colon F \rightarrow D_ {E_1} $ such that $ p_ {E, F}$ takes faces to faces and
$ \chi_E | _F = \chi_ {E_1} \circ p_ {E, F} $;
\item \label{wp3}
for any faces $ F '<F <D_E $ and $ F_1 = p_ {E, F} (F') <D_ {E_1} = p_ {E, F} (F) $, the equality $ p_ {E, F '} = p_ {E_1, F_1} \circ (p_ {E, F} | _ {F '}) $ holds.
\end{enumerate}
Then $ X $ is called a \emph{weakly polyhedral cellular space}.
\end{opr}
By virtue of the above conditions, the definition domains  of the characteristic mappings for all cells of the space $ X $ can be glued using the gluing mappings $ p_ {E, F} $. It is clear that the resulting space is homeomorphic to $ X $.

\begin{pr} \label{Hwp}
The space $ \bar {H} _w $ and its subspace $ \bar {H} _w \smallsetminus H_w $ are weakly polyhedral.
\end{pr}
\begin{proof}
By Corollary \ref {cellsubspace} the space $ \bar {H} _w $ is cellular and by Theorem \ref {Clcell} is framed. For its cell $ E $ and a face $ F <D_E $, put $ p_ {E, F}=\prr_A|_F $ where $ \prr_A $ is the projection from Proposition \ref {prF} if $ F $ is exceptional, and put $ A = \varnothing $ and $ p_ {E, F} = \mathrm {id} $ if  $ F $ is non-exceptional. The conditions of the previous definition are satisfied: (\ref {wp1}) is obvious; (\ref {wp2}) follows from Remark \ref {restr} and the fact that $ \prr_A $ takes  faces to  faces as seen from the proof of Theorem \ref {Clcell}; (\ref {wp3}) follows from the equality $ \prr_B \circ \prr_A = \prr_ {A \cup B} $ for $ B \subset \{1,2, \ldots, a \} \smallsetminus A $. The same is true for $ \bar {H} _w \smallsetminus H_w $.
 \end{proof}
\subsubsection {A triangulation of a weakly polyhedral cellular space.}
\begin{lem}
\label{Cyl}
Let $ D $ be a compact convex set in a Euclidean space, $ \chi \colon D \rightarrow E $ be a continuous surjection onto a Hausdorff space that is a topological embedding of the interior $ \Int D $ of  $ D, $ and $ C = \Cyl \chi | _ {\partial D} $ be the cylinder of the  restriction of the mapping $ \chi $. Then there is a homeomorphism $ f \colon D \cup C \rightarrow E $.
\end{lem}
\begin{proof}
Choose a point $ O \in \Int D $ and consider the homothety $ h_t \colon D = D_1 \rightarrow D_t $ with the center at  $ O $ and with a coefficient $ t $. Considering the identification of $ (x, 0) $ with $ x $ and $ (x, 1) $ with $ \chi (x) $ $ \forall x \in \partial D $ in $ C $, put $ f (z) = \chi (h_ {1/2} (z))$ for $ z \in D $, and $ f (z) = \chi (h _ {(1 + t) / 2} (x))$ for $ z = (x, t) \in C $. It is clear that $ D \cup C $ is compact, $ f $ is a continuous bijection, and therefore a homeomorphism, since $ E $ is Hausdorff.
 \end{proof}
\begin{te}
\label{G}
A weakly polyhedral cellular space $ X $ is triangulable.
\end{te}
\begin{proof}
For each cell $ E $ we choose a point $ a_E \in \Int D_E $, and for each face $ F <D_E $ and the corresponding cell $ E_1 $ a point $ a_F \in \Int F $ with $ p_ {E, F} (a_F) = a_ {E_1} $. The stellar subdivision as a result of starring  at these points gives a triangulation of all polyhedra $ D_E $.
By Lemma \ref {Cyl}, to triangulate the closed cell $ E $, it is sufficient to triangulate the space $ D_E \cup \Cyl \chi_E | _ {\partial D_E} $. We  use the triangulation of the cylinder of a simplicial mapping proposed in \cite [Sec. 4] {C}. It is certain to be determined after one more stellar subdivision of all polytopes $ D_E $.

Due to  conditions (\ref {wp2}) and (\ref {wp3}) of Definition \ref {wp} the  triangulations of cells induce coinciding triangulations at their intersections and therefore define a triangulation of  $ X $.
\end{proof}
\begin{pr} \label{dretr}
Let $ X $ be a weakly polyhedral cellular space with the triangulation indicated in the proof of Theorem \ref {G}, $ Y $ be its cellular subspace, and $ N $ be a (simplicial) open regular neighborhood of this subspace. Then $ X \smallsetminus N $ is a simplicial subspace of $ X $, which is a deformation retract of $ X \smallsetminus Y $.
\end{pr}
\begin{proof}
It follows from the construction of a triangulation of  $ X $ that if all vertices of some simplex lie in $ Y $, then the whole simplex also lies in $ Y $. Therefore, \cite [Lemma 70.1] {M} gives the required statement.
\end{proof}
\subsubsection {Dual partition} 
Throughout this subsection, $ X $ is a triangulated space.

\emph{The dual partition} of  $ X $ is its partition into barycentric stars (see \cite [section 2.2.6.6] {RF} or \cite [section 8.3] {VF}).
The union of stars of dimension no higher than $ m $
is called the \textit {$ m $-dimensional skeleton of the dual partition}.

The next remark follows from the fact that the link of a simplex is homeomorphic to the barycentric link of this simplex (see \cite [section 2.2.6.7] {RF}).
\begin{za}
\label{Skedual}
In $ X $, if for simplices of dimension
$ \dim X-2 $
their links are 
path-connected then
\begin{itemize}
\item
an open star $ s $ from the dual partition of
$ X $
 is a cell if $ \dim s\leq 2 $;

\item
the two-dimensional skeleton $ \mathrm {Sk} _2DX $ of dual partition of
$ X $
is a cellular space.
\end{itemize}

\end{za}
\begin{te}
\label{pi}
$ \pi_1 (X) = \pi_1 (\mathrm {Sk} _2DX) $.
\end{te}
\begin{proof}
Since $ X $ is triangulated, we can assume that the image of homotopy of any loop does not fill the entire $ d $-dimensional star for $ d> 2 $. Therefore, the homotopy can be contracted to $ \mathrm {Sk} _2DX $.
 \end{proof}
\subsection {Dual partition of $ \bar {H} _w $}

\begin{te}
\label{FrSt}
In the space $ \bar {H} _w $, links of simplices of
 dimension
$ \dim _ {\mathbb {R}} \bar {H} _w-2 $ are path-connected.
\end{te}
\begin{proof}
Let $ s \subset \bar {H} _w $ be a simplex of dimension $ \dim _ {\mathbb {R}} \bar {H} _w-2 $ and $ s_1, s_2 \subset \bar {H} _w $ be  any two simplices  of dimension $ \dim _ {\mathbb {R}} \bar {H} _w $ adjacent to $ s $. For the proof, it suffices to 
connect $ s_1$ with $ s_2$ by a chain $ c $ of simplices adjacent to $ s $ such that two successive simplices in $ c $ have a common face of dimension $ \dim _ {\mathbb {R}} \bar {H} _w-1 $ adjacent to $ s $. Thus, the problem is local in the sense that in a small neighborhood of a point $ [f] $ from the open part of $ s $, it suffices to find a path connecting the points $ [f_i] \in s_i, i = 1,2 $. It is clear that the coverings $ f_i $ can be considered  the perturbations of the covering $ f $. Below we  choose a common perturbation domain of  $ f $, used to obtain $ f_1, f_2 $.

Due to Corollary \ref {OpenCell}
in the target sphere of the mapping~$ f $ the following cases are possible.
\begin{enumerate}
\item
There are two mobile meridians, each containing exactly  two simple mobile values.
\item
There is a mobile meridian containing exactly two simple mobile values and a stationary meridian containing a single simple mobile value.
\item
There are two stationary meridians (possibly coinciding), each containing a single simple mobile value.
\item There is a single double mobile value.
\item There is a single stationary value whose perturbation gives two stationary values and one simple mobile value.
\end{enumerate}
In all cases, select disjoint disks on the sphere, each  containing exactly one specified value. These disks define a perturbation domain  of $ f $. In the first four cases, the existence of a path connecting the points $ [f_i] \in s_i, i = 1,2 $ follows from Proposition \ref {pert}, and in the last case, from Proposition \ref {stat}.
 \end{proof}
The following result is immediately obtained from Proposition \ref {dretr}, Remark \ref {Skedual} and Theorems \ref {pi}, \ref {FrSt}.
\begin{sle}
\label{SkedualH}
Let $ Y \subset \bar {H} _w $ be a cellular subspace ( e.g., $ Y = \bar {H} _w \smallsetminus H_w $). Remove from $ \bar {H} _w $ the open stars of all simplices in $ Y $. We obtain a simplicial subspace such that the two-dimensional skeleton of its dual partition is a cellular space. The fundamental group of this skeleton coincides with the fundamental group of  $ \bar {H} _w \smallsetminus Y $.
\qed
\end{sle}
The next sentence is not used further and is presented for completeness.
\begin{pr}
\label{pseudom}
The space $ \bar {H} _w $
is a pseudomanifold, that is, a dimensionally homogeneous unbranched  strongly connected simplicial space \emph{(see \cite [section 8.1] {VF})}.
\end{pr}
\begin{proof}
Let $ d = \dim \bar {H} _w $. The space $ \bar {H} _w $
is dimensionally homogeneous since any cell of dimension $ <d $ obviously lies on the boundary of a $ d $-dimensional cell.
While a cell can lie on the boundary of a single cell, any $ (d-1) $-dimensional simplex is a face of exactly two $ d $-dimensional simplices  by virtue of the above construction of a triangulation of $ \bar {H} _w $. So this space is unbranched. According to Corollary \ref {SkedualH}, the connectedness of  $ \bar {H} _w $ implies the connectedness of the one-dimensional skeleton of the dual partition of this space, which means its strong connection.
\end{proof}
\section {Trigonal Curves}
\label{Hirz}
Let $ S_e $ be a (complex) Hirzebruch surface, i.e.  a rational ruled surface, $ q \colon S _e \rightarrow \Cb \mathrm {P} ^ 1 $ be the corresponding $ \Cb \mathrm {P} ^ 1 $-fibration with exceptional section $ s $, $ s ^ 2 = -e <0 $. The fibers of the bundle $ q $ are called \emph{vertical lines}.

When contracting an exceptional section to a point, the surface $ S_e $ turns into a weighted projective plane $ \mathrm {P} (1,1, e) $ with coordinates $ x_0, x_1, y $, having weights $ 1,1, e $ (see,  e.g., \cite [1.2.3] {Dol}). In what follows, it is convenient to assume that these coordinates are defined on the surface~$ S_e $ itself. A curve $ C \subset S_e $,
given by an
equation
\begin{equation} y ^ 3 + b (x) y + w (x) = 0, \label{1}
\end{equation}
where $ x = (x_0, x_1) $, and $ b $ and $ w $ are homogeneous polynomials of degrees $ 2e $ and $ 3e $, is called a \emph{trigonal curve}. In this case, we admit the cases $ b = 0 $, $ w = 0 $, but not $ b = w = 0 $.
The polynomials $ b $, $ w $ are uniquely determined by the curve $ C $ up to the transformation
\begin{equation} (b, w) \mapsto (t ^ 2b, t ^ 3w), \, t \in \Cb ^ *, \label{2} \end{equation}
therefore the set of all trigonal curves on $ S_e $
is a weighted projective space $ T_e = \mathrm {P} (2, \ldots, 2,3, \ldots, 3) = \mathrm {P} (2_ {2e + 1}, 3_ {3e + 1}) $ of complex  dimension $ 5e + 1 $.

Let $ d = 4b ^ 3 + 27w ^ 2 $ be the discriminant with respect to $ y $ of the equation (\ref {1}).

A generalization of a trigonal curve is a \emph{trinomial curve} on the surface $ S_e $,
given by an
equation
$ y ^ n + b (x) y ^ {n-k} + w (x) = 0 $,
where $ 1 \! \leq \! k \! <\! n $, $ {\texttt {gcd} (k, n) \! = \! 1} $, and $ b $ and $ w $ are homogeneous polynomials of degrees $ ke $ and $ ne $. Let $ d (x) = n ^ nw ^ k - (- 1) ^ n (nk) ^ {nk} k ^ kb ^ n $ for $ n> k $ and $ d (x) = w + b $ for $ n = k = 1 $.
All the results obtained below for trigonal curves are transferred, up to notation, to trinomial curves.

The following lemmas can be proved by a
standard calculation.
\begin{lem}
\label{sing}
The point $ (x, y) \in C $ is singular if and only if  either $ b(x)\neq 0$, $ 2y ^ 3 = w (x) \neq 0 $, and $ x $ is a multiple root of $ d $ or
$ y = 0, $ $ x $ is a root of $ b $ and the multiple root of $ w $. \qed
\end{lem}

\begin{lem}
\label{inf}
Let $ x $ be a common root of $ b, $ $ w $. A point $ (x, y) \in C $ is non-singular if and only if  $ x $ is a simple root of $ w $. In this case, $ (x, y) $ is a point of triple intersection with the vertical line and is called a \emph{vertical inflection point}. \qed
\end{lem}

As in \cite [3.3.1] {DIK}, a non-singular trigonal curve $ C \subset S_e $ is \emph{general} if $ b $, $ w $ have no common roots and all the roots of these polynomials are simple, and \emph{almost general} if only the first of these two conditions is satisfied.

\subsection {The $ j $-invariant of a trigonal curve.}
The function $ j_C \colon \Cb \mathrm {P} ^ 1 \rightarrow \Cb \mathrm {P} ^ 1 $,
$ j_C = \frac {4b ^ 3} {d} = 1- \frac {27w ^ 2} {d} $,
is called the \emph{$ j $-invariant} of a trigonal curve $ C $.

If $ b $, $ w $ have no common roots, then the degree of the mapping $ j_C $ is $ 6e $. In general, this degree $ s = 6e- \deg \texttt {gcd} (b ^ 3, w ^ 2) $ defines a stratification $ \{T_ {e, s} \} $ of the space $ T_e $.
	The non-singular curves of the highest stratum $ T_ {e, 6e} $ form the set of all almost general curves.
	
	\begin{pr}
	\label{Te6eConnect}
	The set $ T_ {e, 6e} $ is connected.
	\end{pr}
	\begin{proof}
	The set $ T_ {e, 6e} $ is determined in the connected space $ T_e $ by the condition that the resultant of  $ b $, $ w $ is nonzero  and therefore is connected.
	\end{proof}
\subsection {Curves with the constant $ j $-invariant.}

Following  \cite [3.1.1] {Degt}, we call the curves $ C \in T_ {e, 0} $, i.e. the curves with the constant $ j $-invariant, \textit {isotrivial}, and with the non-constant one \textit {non-isotrivial}.
\begin{za}
\label{SingTrig}
According to {\rm \cite {DD}}, the space $ T_e = \mathrm {P} (2_ {2e + 1}, 3_ {3e + 1}) $ is singular at a point with vanishing either the first $ 2e + 1 $ coordinates or the last $ 3e + 1 $ ones, i.e.  at a point corresponding to the trigonal curve with $ b = 0 $ or $ w = 0 $. Therefore, the space $ T_e \setminus T_ {e, 0} $ of  non-isotrivial curves is non-singular, i.e. is a manifold.
\end{za}
The next statement is obvious.
\begin{pr}
\label{J_0}
The set ~$ T_ {e, 0} $  is the union of the subspaces ~$ {b = 0}, $ $ w = 0 $ and the submanifold consisting of curves with $ b = \lambda a ^ 2, $ $ w = \mu a ^ 3, $ where $ \lambda, \mu \in \Cb $ and $ a $ is a polynomial of degree $ e $.
\qed
\end{pr}
\subsection {Riemann data of trigonal curve}
\label{RD}

For the hosohedron $ \Theta (j_C) $ of a non-isotrivial trigonal curve $ C $, add to its stationary vertices $ 0 $ (\black--vertex) and $ \infty $ (\cross--vertex) a point $ 1 $ and call it \white--\emph{vertex}. Consider the meridian passing through  $ 1 $  as the prime meridian. Denote by $ [0, \infty] _C $ the subgraph of $ \Theta (j_C) $, consisting of the prime meridian and the vertices lying on it, including the poles. By the
\emph{ graph $ \Gamma (C) $ of trigonal curve} $ C $ we mean the graph $ \Gamma (j_C) \cup j_C ^ {- 1} ([0, \infty] _C) $. In this case, all roots of the polynomial $ w $ are called \white--\emph{vertices} and are added to the stationary vertices of the union, and possible mobile vertices of degree $ 2 $  on the $ 0 $-th monochrome part are removed. Thus, $ \Gamma (C) $ is a $ 0 $-in-graph and a hosograph.

Recall that by definition (see Subsection \ref {Deform}) the index of a point $ v $ on the sphere,  the definition domain of the $ j $-invariant, is equal to $
\mathrm {deg} _vj_C $. Therefore, the index of a mobile and a \white--vertex of $ \Gamma (C) $ equals half of its valency, and of stationary one equals its valency, divided by $ a_ {j_C} $.

It follows from the definition of the $ j $-invariant that for $ C \in T_ {e, 6e} $ the indices of \black--vertices $ \equiv 0 \mod 3 $ and of \white--vertices $ \equiv 0 \mod 2 $.

Let $ \mathrm {cr} (C) $ be
the multiset consisting of all critical values of  $ j_C $, taken with their multiplicities, (i.e.  the function given by $ t \mapsto \sum_ {x \in f ^ {- 1} (t)} (\ind x-1 ) $, where $ t $ is a critical value of $ j_C $, and $ \ind $ is the index of $ x $).
The pair $ (\Gamma (C), \mathrm {cr} (C)) $ is called the \textit {Riemann data} of  $ C $ (cf. \cite [section 1.8] {LZ}).

Let $ \Gamma $ be an $ a$-$ 0 $-in-graph and $ f \colon \Gamma \rightarrow \Theta $ be a morphism onto a $ 0 $-in-hosohedron. Take an interior point of the prime meridian and  call it a \white--vertex.
Add it
to the stationary vertices of $ \Theta $ and call all its preimages in $ \Gamma $ stationary \white--vertices. The graph $ \Gamma $ is called a \white--\emph{graph} (more precisely, an $ a$-\white--\emph{graph}). It is clear that the graph of a trigonal curve  is a \white--graph. We call \black--vertex (respectively, \white--vertex) of a \white--graph \emph{exceptional} (compare \cite [3.1.2] {Degt}) if its index  $ \not \equiv 0 \mod 3 $ (respectively, $ \not \equiv 0 \mod 2 $).

A multiset $ \mathrm {cr} $ of points of the sphere and a \white--graph $ \Gamma $ are called \emph{consistent} if there exists a morphism $ f \colon \Gamma \rightarrow \Theta $, where $ \Theta $ is a hosohedron defined by the multiset $ \mathrm {cr} $, and the multiplicity of each point $ t \in \mathrm {cr} $ is $ \sum_ {x \in f ^ {- 1} (t)} (\ind x-1 ) $.

The following construction recovers a non-isotrivial trigonal curve from its Riemann data.
Take a
natural number $ e $, an $ a$-\white--graph $ \Gamma $ with $ (a + 1) s $ faces, where $ s \leq 6e $, and a multiset $ \mathrm { cr} $ of points of the sphere, consistent with $ \Gamma $. Since $ \Gamma $ is a hosograph, by Theorem \ref {graphf} and \cite [Theorem~3] {ZO} there exists a rational function $ j = b_1 / d_1 = 1-w_1 / d_1 $ of degree~$ s $ with $ \texttt {gcd} (b_1, w_1) = 1 $,  $ \Gamma (j) = \Gamma $ and  the multiset of critical points $ \mathrm {cr} $.
Since $ j $ is not constant, according to \cite [Theorem 3.20, Remark 3.21] {Degt} there exists a polynomial $ g $, whose roots with their multiplicities are exceptional vertices of $ \Gamma $, and there are polynomials $ b $, $ w $ with $ b_1g = 4b ^ 3 $, $ w_1g = 27w ^ 2 $, $ \deg b = 2e $, $ \deg w = 3e $. They define a non-isotrivial trigonal curve $ C \in T_ {e, s} $ such that $ j = j_C $, $ g = \texttt {gcd} (b ^ 3, w ^ 2) $ and $ \Gamma = \Gamma (C) $.

The curve $ C $ is unique up to the action of the group $ G = PGL (2, \Cb) $, and the polynomials $ b $, $ w $ are uniquely determined up to the transformation ~(\ref {2}). Therefore, \textit {the set of Riemann data of trigonal curves can be identified with the quotient space} $ RD_e = (T_e \setminus T_ {e, 0}) / G $. It is clear that its complex dimension is $ 5e-2 $.
 By $ rd \colon T_e \setminus T_ {e, 0} \rightarrow RD_e $ we denote the restriction of the projection $ pr \colon T_e \rightarrow T_e / G $ to the set $ T_e \setminus T_ {e, 0} $ .

\begin{figure} [ht]
\includegraphics {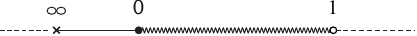} \\
\caption {Coloring of the real projective line.} \label{RP1}
\end{figure}

The following example shows that the space of Riemann data  is not Hausdorff. Take Riemann data with the graph shown in Figure \ref {VertFlex}.a (in this and the following figures, it is assumed that the graph $ \Theta (j) $ is colored as in Figure~\ref {RP1}),
and  tend $ j (x_2) $ to $ j (x_1) = 0 $. In the limit, we get Riemann data both with the graph \ref {VertFlex}.b, and with the graph \ref {VertFlex}.c. Indeed, if for a trigonal curve $ C $ with the graph \ref {VertFlex}.a we assume that the \black--vertex $ x_1 $  equals $ 0 $, the \white--vertex inside $ S_1 $  equals $ t ^ 2 $, and the second \black--vertex is  $ \infty $, then in affine coordinates $ j_C (x) = \frac {4x ^ 3} {4x ^ 3 + 27w_0 (xt ^ 2) ^ 2 (x- \lambda) ^ 2 (x- \mu) ^ 2} $. On the other hand, the $ j $-invariant of the curve $ C_t $,
obtained from $ C $ by the linear fractional transformation $ x \mapsto t ^ 2x / ((1-t) x + t) $ and therefore lying together with $ C $ in the same orbit of the group $ PGL (2, \Cb) $,  equals $
\frac {4x ^ 3 (x-tx + t) ^ 3} {4x ^ 3 (x-tx + t) ^ 3 + 27w_0 (x-1) ^ 2 (\lambda x + o (t)) ^ 2 (\mu x + o (t)) ^ 2} $. Therefore, for $ t \rightarrow 0 $ this orbit approaches curves with the graphs \ref {VertFlex}.b and \ref {VertFlex}.c.
\begin{figure} [h]
\begin{center}
\scalebox {1} {\includegraphics {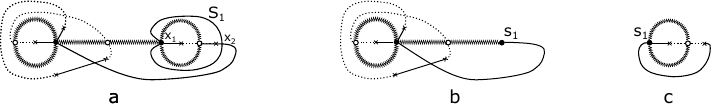}} \\
\end{center}
\caption {Degenerations of the \black--vertex of the graph {\tt a} into the exceptional vertices of the graphs~{\tt b} ~and ~{\tt c}.}
\label{VertFlex}
\end{figure}

\vspace {-3mm}
\subsection {The compactification  of the senior stratum  for the space  of Riemann data of trigonal curves}
\label{CompactRD}
A surface that is a degeneration of a sphere is called a \emph{surface of genus} $ 0 $. For such a surface $ L $ consider a graph whose vertices are the singular points of $ L $ named $ s $\nobreakdash-\emph{vertices} and the components of $ L $ named $ c $-\emph{vertices}, and whose edges connect a component with the lying on it singular points. For a singular surface, this graph is bipartite. We call it the $ cs $-\emph{graph} of $ L $.

We denote by $ F_n $ the space of complex rational functions $ f \colon \Cb \mathrm {P} ^ 1 \rightarrow \Cb \mathrm {P} ^ 1 $ of degree $ n $ and by $ \mathcal {X} _n $ the quotient space $ F_n / PGL (2, \Cb) $. Due to \cite [Theorem 3] {ZO} and Corollary \ref{compact}, $ \bar {\mathcal {X}} _ {n} $ \emph{is a finite cellular space.}  Let 
$ J_ {e} $ be the set of isomorphism classes of $ j $-invariants of curves from~$ T_ {e, 6e} $.
The arguments in Subsection \ref {RD} and \cite [Theorem 3] {ZO} show that $ J_ {e} $ is homeomorphic to $ rd (T_ {e, 6e}) $ and is contained in $ \mathcal {X} _ { 6e} $. The closure $ \bar {J} _ {e} \subset \bar {\mathcal {X}} _ {6e} $ of this set consists of the classes $ [j] $ such that  $ j $ is a branched covering of degree $ 6e $,  the number of \black--vertices of  its graph $ \Gamma (j) $ is  $ \leq 2e $   and the number of \white--vertices is $ \leq 3e $. Therefore, $ \bar {J} _ {e} $ \emph{is a cellular subspace of} $\bar {\mathcal {X}} _ {6e} $ by virtue of Corollary \ref {cellsubspace}, since $ \bar {J} _ {e} = \bar {H} _ {w} $, where $ w $ = (\black -, \white-).

Define a curve $ C $, composed of trigonal curves, whose base is a singular rational curve $ L $. Let $ L $ splits into irreducible components $ L_1, \ldots, L_m $. For each $ i = 1, \ldots, m $, consider the Hirzebruch surface $ S_ {e_i} $ with the base $ L_i $ and a trigonal curves $ C_i \subset S_ {e_i} $, satisfying the following conditions:
\smallskip

\noindent 1) over a common point of  $ L_i $, $ L_l $ there is a common fiber of the surfaces $ S_ {e_i} $, $ S_ {e_l} $ that intersects each of the curves $ C_i $, $ C_l $ at their common triple singular point $ (x_ {il}, 0) $;

\smallskip

\noindent 2) the $ j $-invariants of  $ C_i $, $ C_l $ coincide at the point $ x_ {il} $.\smallskip\\ 
It is easy to check that under these conditions the point $ x_ {il} $ is a common root of  $ \texttt {gcd} (b_i ^ 3, w_i ^ 2) $, $ \texttt {gcd} (b_l ^ 3, w_l ^ 2) $, where $ b_i $, $ w_i $ and $ b_l $, $ w_l $ are the polynomials defining the curves $ C_i $ and $ C_l $.
We call the curve $ C = \bigcup_ {i = 1} ^ m C_i $ \textit {compound (more precisely, $ m $-compound) trigonal curve} with the base $ L $. A trigonal curve with a non-singular base is a $ 1 $-compound trigonal curve.

\smallskip

Gluing the $ j $-invariants of the
curves $ C_i $ gives the $ j $-invariant $ j_C \colon L \rightarrow \Cb \mathrm {P} ^ 1 $ of the curve $ C $. Therefore, the Riemann data of compound trigonal curve are determined.
\emph{The degree of a compound trigonal curve} is the degree of $ j_C $. A perturbation/degeneration of a trigonal curve is a perturbation/degeneration of its $ j $-invariant.

Let $ \bar {T} _ {e, 6e} $ be the set of compound trigonal curves of degree $ 6e $. In view of the above, $ rd (\bar {T} _ {e, 6e}) $ is a compactification of the space $ rd (T_ {e, 6e}) $. We identify $ J_ {e} $ with $ rd (T_ {e, 6e}) $ and $ \bar {J} _ {e} $ with $ rd (\bar {T} _ {e, 6e}) $.

\begin{figure} [ht]
\hspace {15mm} \scalebox {0.7} {\includegraphics {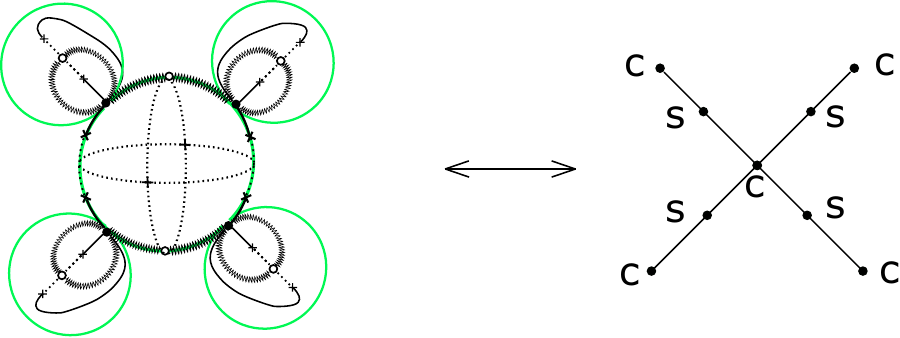}}
\caption{The graph of a compound trigonal curve (left) and the $ cs $\nobreakdash-graph of the base of the curve (right).} 
\label{cs-graph}
\end{figure}

According to Lemmas \ref {sing}, \ref {inf}, the vertical inflection points of a trigonal curve are nonsingular, but a curve with such points does not lie in $ T_ {e, 6e} $. If the other points of the curve are non-singular, it is a component of a compound trigonal curve $ C \in \bar {T} _ {e, 6e} \setminus T_ {e, 6e} $. The $ cs $\nobreakdash-graph of the base~$ L $ of $ C $ is a star \emph{with simplest ends}, i.e.  its end vertices correspond to the components of  $ L $, on which $ j_C $ has degree $ 4 $ (see Figure \ref{cs-graph} with the graph of a compound trigonal curve and the $ cs $\nobreakdash-graph of the base of the curve). All curves in $ \bar {T} _ {e, 6e} \setminus T_ {e, 6e} $ with the $ cs $\nobreakdash-graph not a star with simplest ends we  call \emph{singular compound trigonal curves} and add them to singular curves of the stratum $ T_ {e, 6e} $, obtaining the set $ ST_e $ of singular curves of the space $ \bar {T} _ {e, 6e} $.
\begin{pr}
\label{singspace}
The set $ rd (ST_e) $ is a cellular subspace of $ rd (\bar {T} _ {e, 6e}) $.
\end{pr}
\begin{proof}
A curve $ C \in \bar {T} _ {e, 6e} $ is singular if its discriminant $ d $ has multiple roots  and the $ cs $-graph of its base is not a star with simplest ends.  An open cell in $ rd (ST_e) $ consists of the classes of singular curves. The boundary points of the cell correspond to degenerations of these curves. Such a degeneration does not eliminate multiple roots of the discriminant and does not transform the $ cs $-graph  into a star  with the simplest ends. Therefore the cell boundary also lies in $ rd (ST_e) $.
\end{proof}

\subsection {The dual partition of the space of Riemann data}

\begin{te}
\label{FrStT}
The link of a $ 2 $-dimensional star from the dual partition of the space $ rd (\bar {T} _ {e, 6e}) $ is path-connected.
\end{te}
\begin{proof}
In Theorem \ref{FrSt}, we can take $ rd (\bar{T} _ {e, 6e})$ as $ \bar{H} _w. $
\end{proof}
\vspace{1ex}
From Propositions \ref {Hwp}, \ref {singspace} and Theorem \ref {FrStT}, we obtain the following result.
\begin{sle}
\label{dualsing}
The fundamental group of the space of almost general curves \textup (respectively, of nonsingular curves \textup) is calculated according to the scheme described in Corollary \ref {SkedualH} . \qed
\end{sle}

The authors express their gratitude to the referee for advice and comments that made it possible to eliminate the shortcomings in the original version of the article.

\end{document}